\documentclass[11pt]{article}
\usepackage{geometry}  
\usepackage{graphicx}
\usepackage{amssymb}
\usepackage{epstopdf}
\usepackage{amsmath}
\usepackage{amssymb}
\usepackage{amsthm}
\newtheorem{theorem}{Theorem}
\newtheorem{lemma}[theorem]{Lemma}

\newtheorem{fact}[theorem]{Fact}

\newtheorem{conjecture}[theorem]{Conjecture}

\newtheorem{problem}[theorem]{Problem}

\title{Antipaths in oriented graphs}

\providecommand{\keywords}[1]{{\textit{Keywords:}} #1}

\author{\smallskip Tereza Klimo\v{s}ov\'a\footnote{TK is affiliated to Department of Applied Mathematics, Faculty of Mathematics and Physics, Charles University, Prague, Czech Republic and is supported by the Center for Foundations of Modern Computer Science (Charles Univ. project UNCE/SCI/004) and by GACR grant 22-19073S.  Email: \texttt{tereza@kam.mff.cuni.cz}} \
and\ Maya Stein\footnote{MS is affiliated to the Department of Mathematical Engineering and the Center for Mathematical Modeling of the University of Chile, and acknowledges support by FONDECYT Regular Grant 1221905, by FAPESP-CONICYT Investigaci\'on Conjunta Grant C\'odigo 2019/13364-7, and by ANID Basal Grant CMM FB210005. Email: \texttt{mstein@dim.uchile.cl}}}

\date{}                                           

\begin{document}
\maketitle
\begin{abstract}
We show that for any natural number $k\ge 1$, any oriented graph $D$ of minimum semidegree at least $(3k-2)/4$ contains an antidirected path of length $k$. 

In fact, a slightly weaker condition on the semidegree sequence of $D$ suffices, and as a consequence, we confirm a weakened antidirected path version of a conjecture of Addario-Berry, Havet, Linhares Sales, Thomass\'e and Reed. 
\end{abstract}


\section{Introduction}
In undirected graphs, a large minimum degree is very helpful for finding long paths. For instance, if we wish to ensure that an $n$-vertex   graph $G$ contains a path of length $k$ (i.e.,~with $k$ edges), a simple greedy embedding argument shows that it is enough to assume $G$ has minimum degree at least $k$. Although in general, this bound on the minimum degree is best possible, there is room for improvement if $k$ is large compared to $n$, or if  $G$ is assumed to be connected. Dirac~\cite{dirac} showed that a minimum degree of at least $n/2$ is enough to find a Hamilton cycle  in an $n$-vertex graph, and with a minimum degree exceeding $(n-2)/2$ we can find a Hamilton path, i.e.,~a  path of length $n-1$. With a similar proof, one can show that every connected graph on at least $k+1$ vertices that has minimum degree strictly greater than  $(k-1)/2$ contains a path of length $k$.

It would be interesting to find extensions of these results to digraphs. We will focus on oriented graphs here (see Section~\ref{sec:final} for some remarks on the general digraph case).
 We have to decide which parameter will play the role of the minimum degree, and the widespread notion of the {\it minimum semidegree} $\delta^0(D)$, which is defined as the minimum over all out-and in-degrees of all vertices of the oriented graph $D$,  seems a natural choice. 
In the same way as in the undirected case, we can use a greedy embedding strategy  to see that  any oriented graph $D$ with $\delta^0(D)\ge k$ must contain each orientation of the $k$-edge path. And as before, it seems reasonable to ask whether  this bound can  be lowered if the underlying graph $G$ of~$D$ (i.e.,~the graph we obtain by omitting directions) has a sufficiently large connected component. Note that the condition $\delta^0(D)\ge k/2$ alone already implies that $G$ has a connected component with at least $k+1$ vertices.

There are many results for Hamilton cycles  in oriented graphs. As a Hamilton cycle of an $n$-vertex oriented graph $D$ contains a path of length $n-1$, these results shed some light on our problem. 
In 1960, Ghoulia-Houri~\cite{g-h} proved a minimum semidegree of $n/2$  guarantees a directed Hamilton cycle in any $n$-vertex digraph $D$, and 
in the 1970's, 
Thomassen~\cite{thomassen}  asked for an analogous result for oriented graph, with a weaker condition on the minimum semidegree. H\"aggkvist~\cite{h} conjectured this  to be $\delta^0(D)\ge (3n-4)/8$, which he showed to be best possible, and which, after previous results in~\cite{ht, kko},  was confirmed by Keevash, K\"uhn and Osthus in~\cite{keeko} for all large oriented graphs. 
H\"aggkvist and Thomassen~\cite{ht2} conjectured that for all $\alpha >0$, all sufficiently large oriented graphs $D$ with $\delta^0(D)\ge (3/8 +\alpha)n$ contain any orientation of a Hamilton cycle, and this was confirmed by Kelly~\cite{kelly}. 
In particular, it follows that for $3n/4+o(n)\le k<n$, every oriented graph on~$n$ vertices and of minimum semidegree at least $k/2$ contains each orientation of the $k$-edge path.

A corresponding result for oriented paths of  length  below $3n/4$ is still missing, except in the case of directed paths: Jackson~\cite{jackson}  showed in 1981 that for every $\ell\in\mathbb N$,  every oriented graph $D$ with $\delta^0(D)\ge\ell$  contains the directed  path on $2\ell$ edges.
In~\cite{trees-survey}, the second author suggested that something similar might be true for all orientations of the $k$-edge path.
\begin{conjecture}$\!\!${\rm\bf\cite{trees-survey}}
\label{conj}\label{q1ori}
For each $k\in\mathbb N$, every oriented graph $D$ with $\delta^0(D)> k/2$  contains each orientation of the  path of length $k$.
\end{conjecture}

 Conjecture~\ref{q1ori}  is sharp in the following sense. The bound on the minimum semidegree could not be lower than $k/2$, as one can see by considering the disjoint union of regular  tournaments on $k$ vertices, if~$k$ is odd.
For {\em antidirected paths} (i.e.,~oriented paths that alternate edge directions), one can also consider the blow-up of a directed cycle of  length~$\ell$, where each vertex $v$ of $C_\ell$ is replaced by an independent set $S_v$ of size $k/2$, and $S_v, S_w$ span a complete bipartite graph whenever $vw\in E(C_{\ell})$. Any largest antidirected path in this graph has length $k-1$, and the minimum semidegree of the graph is $k/2$.
 
As noted above, Conjecture~\ref{q1ori} is true for $n$-vertex oriented graphs and $k\ge 3n/4+o(n)$ by the results of~\cite{kelly}, and it is also true  for directed paths~\cite{jackson}. It has been verified for all $k\le 5$~\cite{cris}.
Z\'arate-Guer\'en and the second author showed in~\cite{camila} that an approximate version of Conjecture~\ref{q1ori} holds  for
 antidirected paths in large oriented graphs~$D$, if $k$ is linear in $n$.

We will focus here on a variant of Conjecture~\ref{q1ori} for
 antidirected paths. 
We show, with a much easier proof than the one from~\cite{camila}, and for any $k\in\mathbb N^+$, that every oriented graph~$D$ with $\delta^0(D)\ge (3k-2)/4$ contains each antidirected path.\footnote{Note that if $k$ is odd, there is only one antidirected path of length $k$ (unless we specify a starting vertex). If $k$ is even, there are two distinct antidirected paths.}
Actually, we will prove a slightly stronger statement. We define the {\it minimum pseudo-semidegree} $\bar\delta^0(D)$ of a digraph~$D$  as follows: $\bar\delta^0(D)=0$ if $D$ has no edges, and otherwise $\bar\delta^0(D)$ is the maximum $d\in\mathbb N$ such that each vertex in $V(D)$ has out-degree either $0$ or $\ge d$, and in-degree either $0$ or $\ge d$. Clearly $\bar\delta^0(D)\ge \delta^0(D)$ for each digraph $D$. Our main result is the following.

\begin{theorem}\label{antii}
Let $k\in\mathbb N$ with $k\ge 3$ and let $D$ be an oriented graph with $\bar \delta^0(D)\ge (3k-2)/4$. Then~$D$ contains each antidirected path of length $k$.
\end{theorem}

Note that the case $k=2$ needs to be excluded from our theorem, because the bound $\bar \delta^0(D)\ge (6-2)/4=1$ is below the bound from Conjecture~\ref{q1ori} and not sufficient to guarantee an antipath of length two (as $D$ could be a directed cycle).

In a similar vein as Conjecture~\ref{q1ori},
Addario-Berry, Havet, Linhares Sales, Thomass\'e and Reed conjectured the following in 2013.
\begin{conjecture}[Addario-Berry et al.~\cite{ABHLSTR}]
\label{antiES}
Every  digraph $D$ with more than $(k-1)|V(D)|$ edges contains each  antidirected tree with $k+1$ vertices. 
 \end{conjecture}
 For symmetric digraphs, this conjecture is equivalent to the Erd\H os-S\'os conjecture; and in oriented graphs, Conjecture~\ref{antiES} implies Burr's conjecture for antidirected trees (see \cite{ABHLSTR} for details). Conjecture~\ref{antiES} is proved in \cite{ABHLSTR} for trees of diameter at most 3, and an approximate version for large balanced antidirected trees in dense oriented graphs is given in~\cite{camila}. 
 
 It is also shown in \cite[Theorem 17]{ABHLSTR}  that every  digraph $D$ with more than $4(m-1)|V(D)|$ edges contains each  antidirected tree whose largest partition class has at most $m$ vertices. This implies that every  digraph $D$ with more than $4(\lceil(k+1)/2\rceil-1)|V(D)|$ (that is, roughly $2k|V(D)|$) edges contains each  antidirected $k$-edge path.
We  improve this bound to roughly $3k|V(D)|/2$. 

 \begin{theorem} \label{corol}\label{antiii}
For each $k\in \mathbb N^+$, every oriented graph $D$ with more than $(3k-4)|V(D)|/2$ edges contains each  antidirected path of length $k$. 
 \end{theorem}





\section{Notation}
A {\it digraph} has directed edges, at most one for each direction between each pair of vertices~$u,v$. For brevity we write {\it edge} instead of {\it directed edge}, and let $uv$ denote an edge going from vertex $u$ to vertex $v$. For the endvertices of such an edge, we say that $v$ is an {\it out-neighbour} of $u$, and $u$ is an {\it in-neighbour} of $v$.
We write $d^-(v)$ and $d^+(v)$ for the  {\it in-degree} and the {\it out-degree} of  vertex $v$: this is the number of out-neighbours, or in-neighbours of $v$, respectively. As already mentiond in the introduction, the {\it minimum semidegree} of a digraph $D$ is $\delta^0(D)=\min\{d^-(v), d^+(v):v\in V(D)\}$, and the {\it minimum pseudo-semidegree} $\bar \delta^0(D)$ of a digraph $D$ is the minimum of $\min\{d^-(v):v\in V(D), d^-(v)>0\}$ and $\min\{d^+(v):v\in V(D), d^+(v)>0\}$, unless $D$ has no edges, in which case $\bar \delta^0(D)=0$.

In an {\it oriented  graph},  for each pair of vertices $u,v$, at most one of the edges $uv$, $vu$ is present.  We say an oriented path or cycle has {\it length} $k$ if it has $k$ edges. 
An {\it antidirected path (antidirected cycle, antidirected tree)} is an oriented path (cycle, tree) where every vertex has either out-degree $0$ or in-degree $0$. We also write {\it antipath (anticycle, antitree)} for short. Note that each anticycle has even length. In particular, any anticycle in an oriented graph has length at least $4$.

\section{Proof of Theorem~\ref{antii}}

We show Theorem~\ref{antii} by combining three auxiliary lemmas, namely Lemmas~\ref{lem1},~\ref{lem2} and~\ref{lem3}, which are stated and proved below. 
The proofs of these lemmas make use of different variants of a well-known argument that appears in the proof of Dirac's theorem. For convenience, we state this tool now, as Fact~\ref{fact},  and include its short proof for completeness. 

Given a set $F$ of edges in an undirected graph, we write $d_F(v,S)$ for the number of edges  between a vertex $v$ and a set $S$ that belong to $F$.
\begin{fact}\label{fact}
Let $m\in \mathbb N^+$, let $1\le \ell\le m$ and let~$G$ be a graph. Let $X,Y\subseteq V(G)$, with $X=\{x_0, x_1, \ldots, x_{m-1}\}$ and $Y=\{y_1, y_2, \ldots, y_m\}$, and let $F_0, F_m\subseteq E(G)$. If  $d_{F_0}(x_0, Y)+d_{F_m}(y_m, X)\ge m+\ell$, then there is an index $i$ with $\ell\le i\le m$ such that $x_0y_i\in F_0$ and $x_{i-\ell}y_m \in F_m$.
\end{fact}
\begin{proof}
Otherwise, $1_{x_0y_i\in F_0}+1_{x_{i-\ell}y_m \in F_m}\le 1$ for each $i=\ell,\ldots,m$, and therefore, $$m+\ell\le d_{F_0}(x_0, Y)+d_{F_m}(y_m, X) \le 2(\ell -1)+ \sum_{i=\ell}^m (1_{x_0y_i\in F_0}+1_{x_{i-\ell}y_m \in F_m})\le m+\ell-1,$$
a contradiction.
\end{proof}

The usual application of this argument in the proof of Dirac's theorem is setting $\ell=1$, and, given a maximum length path $P=x_0x_1\ldots x_m$, setting $y_i:=x_i$ for $i=1,\ldots,m$, while letting $F_0$, $F_m$ be the set of edges going from $x_0$ or $x_m$, respectively, to other vertices of $P$. Then the two edges $x_0x_i$ and $x_{i-1}x_m$ given by Fact~\ref{fact} are used to find a cycle in $V(P)$.

We are ready for the first auxiliary lemma.

\begin{lemma}\label{lem1}
Let $k\in\mathbb N$ and let $D$ be an oriented graph of minimum pseudo-semidegree $\bar\delta^0(D)\ge k/2$. Let $P=v_0v_1\ldots v_m$ be a longest antipath in $D$. If $m<k$ then $m$ is odd.
\end{lemma}
\begin{proof}
Assume otherwise, that is, suppose $m<k$ and $m$ is  even. We may assume that $m\neq 0$ (because $m=0$ means that the longest antipath is trivial, implying that $\bar\delta^0(D)=0$ and thus $k=0$). Then   $m$ fulfils $3\le m+1\le k$.  By symmetry, we may assume that 
\begin{equation}\label{111}
v_0v_1, v_mv_{m-1}\in E(D),
\end{equation}
 that is, the first edge of $P$ is directed towards $v_1$, and that the last edge of $P$ is directed towards $v_{m-1}$. 

Note that by the maximality of $P$, all out-neighbours of $v_0$ lie on $P$, and the same is true for $v_m$. By~\eqref{111}, each of $v_0$, $v_m$ has at least one out-neighbour, and therefore, by our assumption on the minimum pseudo-semidegree of~$D$,  each has at least $k/2\ge (m+1)/2$ out-neighbours. 

Let $G$ be the underlying graph of~$D$, and let $F_0, F_m\subseteq E(G)$ be the sets of all edges of~$G$ corresponding to edges of $D$ that are leaving $v_0$ or leaving  $v_m$, respectively. Then we can calculate $$d_{F_0}(v_0, V(P-v_0))+ d_{F_m}(v_m, V(P-v_m))\ge (m+1)/2+(m+1)/2=m+1.$$
 Let $x_i=y_i=v_i$ for $i=1,\ldots,m$. We now use Fact~\ref{fact} with $\ell=1$ in  $G$ to see that there is an index $i\in\{1,\ldots,m\}$ such that $v_0v_i\in F_0$ and $v_{i-1}v_m\in F_m$. So $v_0v_i$ and $v_mv_{i-1}$  are edges (in these directions) in $D$. 
 
 We may assume that $v_iv_{i-1}$ is an edge (and thus $i\neq 1$), as otherwise we could reverse~$P$, interchanging the roles of $v_0$ and $v_m$, and of $v_i$ and $v_{i-1}$, and thus obtain the desired direction. So  also  $v_iv_{i+1}$ is an edge (unless $i=m$).

Consider the two antipaths $$P'=v_iv_{i+1}\ldots v_mv_{i-1}v_{i-2}\ldots v_1v_0$$ and 
 \begin{equation}
    P''=
    \begin{cases}
      v_iv_0v_1\ldots v_{i-1}v_{m}v_{m-1}\ldots v_{i+1} & \text{if}\ i\neq m \\
    v_iv_0v_1\ldots v_{i-1} & \text{if}\ i= m.
    \end{cases}
  \end{equation}
  Both antipaths $P'$ and $P''$ have the same vertex set as $P$ (and hence maximum length) and start with vertex $v_i$. In $P'$, the first edge is directed away from $v_i$, while in $P''$, the first edge is directed towards $v_i$. Therefore, and by the maximality of the antipath  $P$, all in-neighbours and all out-neighbours of vertex~$v_i$ lie on $P$. 
  
  Also note that $v_i$ has both in- and out-neighbours (namely, $v_0$ and $v_{i-1}$). So, since by assumption $\bar\delta^0(D)\ge k/2$, we know that  $v_i$ has at least $k/2$ in-neighbours and at least $k/2$ out-neighbours, and thus, $k=k/2+k/2\le |V(P-v_i)|=m$, in contradiction to our assumption that $m<k$.
\end{proof}

Our second auxiliary lemma turns the maximum length antipath into an anticycle on the same number of vertices.
For this lemma (and only for this lemma) we need a minimum pseudo-semidegree of $\bar\delta^0(D)\ge (3k-2)/4$.

\begin{lemma}\label{lem2}
Let $k\in\mathbb N$, let $D$ be an oriented graph of minimum pseudo-semidegree $\bar\delta^0(D)\ge (3k-2)/4$, and let $m$ be the maximum length of an antipath in $D$. If $1<m<k$, then~$D$ contains an anticycle of length $m+1$.
\end{lemma}
\begin{proof}
Let $P=v_0v_1v_2\ldots v_m$ be an antipath of maximum length in $D$. 
Assume $1<m<k$. In particular, $(3k-2)/4\ge k/2$ and 
Lemma~\ref{lem1} implies that $m$ is odd. 

By symmetry, and since $m$ is odd, we may assume that $v_0v_1, v_{m-1}v_m\in E(D)$, that is, the first edge is directed towards $v_1$, and the last edge is directed towards $v_{m}$. 
Observe that all edges on $P$ are directed from their endvertex of even index towards their endvertex of odd index. Also observe that by maximality of $P$,  all out-neighbours of $v_0$ and all in-neighbours of $v_m$ lie on $P$.

Let $G$ be the underlying graph of $D$. 
We set $x_i:=v_{2i}$ and $y_{i+1}:=v_{2i+1}$ for all $i=0,\ldots,m'-1$, where $m'=(m+1)/2$. Let $X=\{x_0, x_1, \ldots, x_{m'-1}\}$ and $Y=\{y_1, y_2, \ldots, y_{m'}\}$. Note that $X\cup Y$ is a partition of $V(P)$.

Next, define $F_0$ as the set of all edges of $G$ that correspond to an edge of $D$ leaving $v_0$ and ending at a vertex from $Y$. Further,  let  $F_m$ be the set of all edges of $G$ corresponding to edges of $D$ that start at a vertex from $X$ and end at $v_m$. As by assumption $k\ge m+1$, we know that
\begin{align*}
d_{F_0}(v_0, Y)+d_{F_m}(v_m, X) & = d^+(v_0)- |X\setminus\{v_0\}| + d^-(v_m)- |Y\setminus\{v_m\}|  \\ & \ge (3k-2)/2-(|V(P)|-2)\\ & \ge (3m+1)/2-m+1\\ & \ge m'+1.
\end{align*}

Now, applying Fact~\ref{fact} with $\ell=1$  in  $G$  we find an index $i\in\{1,\ldots,m'\}$ such that  $x_0y_{i}=v_0v_{2i-1}\in F_0$ and $x_{i-1}y_{m'}=v_{2i-2}v_m\in F_m$. So, $$v_0v_1v_2\ldots v_{2i-3}v_{2i-2}v_mv_{m-1}\ldots v_{2i}v_{2i-1}v_0$$ is an anticycle of length $m+1$, which is as desired.
\end{proof}

Our last auxiliary lemma uses the anticycle found in the previous lemma, and turns it into an antipath on more vertices.

\begin{lemma}\label{lem3}
Let  $k\in\mathbb N$ and let $D$ be an oriented graph of minimum pseudo-semidegree $\bar\delta^0(D)> k/2$, and let 
$C$ be an anticycle of length $m+1$  in $D$. If $m<k$, then $D$ has an antipath of length $m+1$.
\end{lemma}
\begin{proof}
Let $C=v_0v_1v_2\ldots v_{m}v_0$. By symmetry, we may assume that $v_0v_1, v_{0}v_{m}\in E(D)$.  
Observe that we may assume the following for all $i=0,1,\ldots , m$:
\begin{equation}\label{neighbeven}
\text{If $i$ is even, then all out-neighbours of $v_i$ lie on $C$; and}
\end{equation}
\begin{equation}\label{neighbodd}
\text{if $i$ is odd, then all in-neighbours of $v_i$ lie on $C$,}
\end{equation}
as any such out- or in-neighbour could be added to $C$ to obtain an antipath of length $m+1$ in $D$, and then we would be done.

In particular, \eqref{neighbeven} and~\eqref{neighbodd} imply that all out-neighbours of $v_0$ and all in-neighbours of~$v_{m}$ lie on the anticycle $C$.
Set $x_i=y_i=v_i$ for $i=0,\ldots, m$, and consider the  underlying graph $G$ of $D$. Let $F_0$ comprise of all edges of $G$ corresponding to edges of $D$ that start at $v_0$ and end on $C$. Let $F_{m}$ contain all edges of $G$ corresponding to edges of $D$ that start on $C$ and end at $v_{m}$. Since $\bar\delta^0(D)> k/2$, and $k\ge m+1$ by assumption, we have that $$d_{F_0}(v_0, \{v_1\ldots, v_m\})+d_{F_m}(v_m, \{v_0\ldots, v_{m-1}\})\ge k+1\ge  m+2.$$

Now, we use Fact~\ref{fact} with $\ell=2$  to see that there is an index $i\in\{2,\ldots,m\}$ such that $v_0v_{i}$ and $v_{i-2}v_{m}$  are edges of $D$, in these directions.

Let us first  assume that $i$ is even, that is, $v_{i-2}v_{i-1}, v_{i}v_{i-1}\in E(D)$. 
Since $v_i$ has both in-and out-neighbours in $D$ (for example $v_0$ and $v_{i-1}$), and because of our assumption on the minimum pseudo-semidegree, we know that $v_i$ has at least $k/2$ in-neighbours and $k/2$ out-neighbours. By \eqref{neighbeven}, all out-neighbours of $v_i$ belong to $V(C)$. So at most $|V(C-v_i)|-k/2=m-k/2\le k/2-1$ vertices of 
 $C$ are in-neighbours of $v_i$, which means that $v_{i}$ has an in-neighbour $x\in V(D)\setminus V(C)$. 

Since $\bar\delta^0(D)> k/2$, and since $x$ has an out-neighbour (namely $v_i$), we know that $x$ has at least $(k+1)/2>(m+1)/2$ out-neighbours in $D$. These cannot all lie in $V(C)$, as otherwise one of them would be a vertex $v_j$ with $j$ odd, a contradiction to~\eqref{neighbodd}. Thus vertex~$x$ has an out-neighbour~$y$ that does not lie on $C$.

Consider the antipath $$P=yxv_{i}v_0v_1\ldots v_{i-3}v_{i-2}v_mv_{m-1}v_{m-2}\ldots v_{i+1}.$$ 
As $x,y\in V(P)\setminus V(C)$, and $V(C)\setminus V(P)=\{v_{i-1}\}$, the antipath $P$ has length $m+1$, which is as desired.

If $i$ is odd, we can find, in a similar way as above, vertices $x$ and $y$ such that $$yxv_{i-2} v_mv_{m-1}\ldots v_iv_0v_1\ldots v_{i-3}$$ is an antipath of length $m+1$. This finishes the proof.
\end{proof}

Now we are ready to prove Theorem~\ref{antii}.

\begin{proof}[Proof of Theorem~\ref{antii}]
Let $m$ be the length of a longest antipath $P$ of $D$. It is easy to see that $m\ge 2$. 

First assume $m<k$. Then by Lemma~\ref{lem2}, $D$ has an anticycle of length $m+1$, and therefore, by Lemma~\ref{lem3} (which can be used since $(3k-2)/4>k/2$ for $k\ge 3$), $D$ has an antipath of length $m+1$, a contradiction to the choice of $m$.

So $m\ge k$. Let $P'$ be an antipath of lenght $k$.
 If $m>k$, or if $m=k$ and $k$ is odd, then~$P$  contains $P'$ (possibly reverting $P$). So we can assume $m=k$ and $k$ is even. 
 Now,  we can apply Lemma~\ref{lem1} with $k'=3k/2-1$, because  $\delta^0(D)\ge (3k-2)/4$ implies that $\delta^0(D)\ge k'/2$, and furthermore, $m=k<k'$ since $k>2$. So $m$ is odd, a contradiction.
\end{proof}

\section{Proof of Theorem~\ref{antiii}}

We will start by proving a  lemma that allows us  to rewrite the condition on the edge density as a condition on  the minimum pseudo-semidegree. This lemma also appears in~\cite{camila}, but for completeness, we include its short proof here.
\begin{lemma}\label{split}
Let $\ell\in\mathbb N$. 
If a digraph $D$ has more than $\ell |V(D)|$ edges, then it contains a digraph $D'$ with $\bar\delta^0(D')\ge (\ell+1)/2$.
\end{lemma}
\begin{proof}
Note that the vertices of $D$ have, on average, in-degree greater than  $\ell$ and out-degree  greater than  $\ell$. Consider the following folklore construction of an auxiliary bipartite graph~$B$ associated to $D$: first,  divide each vertex $v\in V(D)$ into two vertices $v_{in}$ and $v_{out}$, letting $v_{in}$ be adjacent to all edges ending at $v$, and letting $v_{out}$ be adjacent to all edges starting at $v$; second, omit all directions on edges. 

Then the average degree of $B$ is greater than   $\ell$, and a standard argument shows that~$B$ has a non-empty subgraph $B'$ of minimum degree exceeding $\ell/2$ (for this, it suffices to successively delete vertices of degree $\le\ell/2$ and to calculate that we have not deleted the entire graph). Translating $B'$ back to the digraph setting, we see that $D$ has a subdigraph~$D'$ with minimum pseudo-semidegree exceeding $\ell/2$.  
 \end{proof}

Now we are ready to prove Theorem~\ref{antiii}.
\begin{proof}[Proof of Theorem~\ref{antiii}]
Use Lemma~\ref{split} to find a subdigraph $D'$ of $D$ with $\bar\delta^0(D')\ge ((3k-4)/2+1)/2=(3k-2)/4$. As a subdigraph of $D$, also $D'$ is an oriented graph. So
 Theorem~\ref{antii} can be applied to find each antipath of length $k$.
 \end{proof}

%
%
%

\section{Final remarks and open problems}\label{sec:final}

 \paragraph{A lower bound in  Theorem~\ref{antii}.}
We  believe the bound in  Theorem~\ref{antii} is not best possible. We think the lower bound $\delta^0(D)>k/2$ from Conjecture~\ref{q1ori} should be closer to the truth, although we have not been able to improve our result in that direction.
We remark that if one could improve the bound from Lemma~\ref{lem2}, then,  following all steps of our proof, one would automatically obtain an improved bound for Theorem~\ref{antii}.

 \paragraph{Other orientations of the path.}
 In Conjecture~\ref{q1ori}, all possible orientations of the $k$-edge path are considered. As a weakening of the conjecture, we could ask the following.
\begin{problem}\label{problem1}
Does Theorem~\ref{antii}, with $\bar\delta^0$ replaced by $\delta^0$, hold for other types of oriented paths? 
\end{problem}
That is, we ask whether  it is true that  every oriented graph $D$ with $\delta^0(D)> (3k-2)/4$  contains each orientation of the  path of length $k$. (Replacing $\bar\delta^0$  by $\delta^0$ is neccessary because large antidirected complete bipartite graphs only contain antidirected subgraphs, as has been observed in~\cite{ABHLSTR, trees-survey, camila}.)
The fact that Problem~\ref{problem1} holds for  the extreme opposites of possible orientations of  paths -- antipaths and directed paths~\cite{jackson} -- may induce some hope that  Problem~\ref{problem1} is true, whether or not Conjecture~\ref{q1ori}  holds.

 \paragraph{Antitrees.}
It  seems natural to replace anti-paths with anti-trees in Conjecture~\ref{q1ori}. Together with Z\'arate-Guer\'en the second author shows in~\cite{camila} the following result, where an anti-tree is balanced if it has as many vertices of out-degree 0 as vertices of in-degree 0.
\begin{theorem}$\!\!${\rm\bf\cite{camila}}\label{camila}
 For all $\epsilon, c>0$ there is an $n_0$ such that for all $n\ge n_0$ and all $k\ge \epsilon n$, every oriented graph 
$D$ with $\delta^0(D)> (1/2+\epsilon)k$ contains each balanced anti-tree $T$ with $k$ edges and with  maximum degree at most $c\log n$.
\end{theorem}
The proof of Theorem~\ref{camila} uses digraph regularity, but it might be possible to find a simpler proof and/or drop either the approximation or the additional condition on the balancedness and the maximum degree of $T$ if we only look for specific anti-trees.
\begin{problem}
Does every oriented graph 
$D$ with $\delta^0(D)> k/2$ (or $\delta^0(D)> (1/2+\epsilon)k$, or $\delta^0(D)>3k/4$) contain each anti-tree $T$ with $k$ edges, if we add some additional restriction on $T$ (e.g. $T$ is a caterpillar, spider, has small diameter,...)?
\end{problem}

 \paragraph{Digraphs.}
Analogous questions can be asked for digraphs. Observe that as for oriented graphs, a greedy embedding argument gives that $\delta^0(D)\ge k$ is enough to guarantee a copy of any oriented $k$-edge tree $T$ in a digraph $D$. So it seems natural to ask whether this bound can be lowered. 
However, in contrast to the situation in oriented graphs, it will now be necessary to add an additional condition on the order of the largest connected component of the underlying graph of $D$, as the minimum semidegree condition alone is not sufficient to ensure that $|V(D)|\ge |V(T)|$  (since $D$ could be the union of complete digraphs of order $\delta^0(D)+1$). 
For instance, one could ask whether every digraph $D$ with $\delta^0(D)> k/2$ (or some other bound) having a component of size at least $k+1$ contains each oriented $k$-edge path. For this question and  further comments see~\cite{trees-survey}.

\section*{Acknowledgment}
This research was initiated at the MATRIX program ``Structural Graph Theory Down\-under'' in November 2019. The authors are grateful to the organizers and MATRIX's hospitality.

\end{document}